\documentclass{article}

\usepackage[dvipdfmx]{graphicx,xcolor}
\usepackage[fleqn]{amsmath}
\usepackage{newtxtext}
\usepackage[varg]{newtxmath}

\def \F {{\mathbb F}}

\newtheorem{theorem}{Theorem}

\newtheorem{corollary}{Corollary}
\newtheorem{conjecture}{Conjecture}

\title{Pseudorandom binary sequences:
quality measures and number-theoretic constructions}

\author{Arne Winterhof\\
 Johann Radon Institute for Computational and Applied Mathematics\\
 Austrian Academy of Sciences
}

\begin{document}

\maketitle

\begin{abstract}
In this survey we summarize properties of pseudorandomness and non-randomness of some number-theoretic sequences and present results
on their behaviour under the following measures of pseudorandomness: balance, linear
complexity, correlation measure of order $k$, expansion complexity and $2$-adic complexity.
The number-theoretic sequences are the Legendre sequence and the two-prime generator, the Thue-Morse sequence and its sub-sequence along squares, and the prime omega sequences for integers and polynomials.
\end{abstract}




\section{Introduction}

Let 
$${\cal S}=(s_n)_{n=0}^\infty, \quad s_n \in \F_2=\{0,1\}, \quad n=0,1,\ldots,$$ 
be a binary sequence. We call it {\em pseudorandom} if it is deterministically generated but cannot be distinguished from a {\em truly random} sequence.
Pseudorandom sequences are crucial for cryptographic applications such as stream ciphers, see for example \cite{cudire}.  

\subsection{Measures of pseudorandomness}
 There are several measures of pseudorandomness which can be used to detect cryptographically weak sequences including 
 \begin{itemize}
  \item balance,
  \item linear complexity,
  \item maximum-order complexity,
  \item correlation measure of order $k$,
  \item expansion complexity
  \item and $2$-adic complexity.
 \end{itemize}
 These measures are partly not independent and partly complement each other. We will discuss some of their relations.
\subsection{Pseudorandom sequences} 
 We summarize
 results on these measures for the following number-theoretic sequences, 
 \begin{itemize}
 \item the Legendre sequence and the two-prime generator,
 \item the  Thue-Morse sequence and its sub-sequence along squares,
 \item the prime omega sequence modulo $2$ for integers and for polynomials.
 \end{itemize}
 
Each section will focus on one of the above measures of pseudorandomness.

It turns out that 
\begin{itemize}
    \item the Legendre sequence has no obvious flaw (if the period is long enough),
    \item the two-prime generator suffers a large correlation measure of order $4$ and is not pseudorandom,
    \item the Thue-Morse sequence has an undesirable deviation from the expected value $N/2$ of the linear complexity, a large correlation measure of order $2$ and a small expansion complexity,
    and is not suitable in cryptography,
    \item the Thue-Morse sequence along squares seems to be an attractive candidate for cryptography,
    \item the $N$th linear complexity of the prime omega sequence for integers seems to be too regular,
    \item there is no obvious deficiency of the omega sequence for polynomials.
\end{itemize}

For earlier surveys on measures of pseudorandomness see \cite{gy,mewi22,hff,ni03,towi,wi10}.

\section{Balance and definitions of the sequences}
\subsection{Definition of balance and its expected value}
The $N$th {\em balance} $B({\cal S},N)$ of a binary sequence~${\cal S}=(s_n)_{n=0}^\infty$ is
$$B({\cal S},N)=\left|\sum_{n=0}^{N-1}(-1)^{s_n}\right|.$$
The balance of a sequence which is not distinguishable from a random sequence should be of order of magnitude $N^{1/2}$,
see Alon et al.\ \cite[Lemma~12]{al},
or at least
$$N^{o(1)}\ll B({\cal S},N)=o(N).$$
Here we use the notation 
$$f(N)=O(g(N))\Longleftrightarrow  |f(N)|\le cg(N)$$
for some absolute constant $c>0$,
$$f(N)\ll g(N) \Longleftrightarrow f(N)=O(g(N))$$
and 
$$f(N)=o(g(N)) \Longleftrightarrow \lim_{N\rightarrow \infty} \frac{f(N)}{g(N)}=0.$$

Now we give a list of some number-theoretic sequences with desirable balance.

\subsection{Legendre sequence}

For a prime $p>2$ the {\em Legendre sequence} ${\cal L}_p=(\ell_n)_{n=0}^\infty$ is the $p$-periodic sequence defined by
\begin{equation}\label{legdef}\ell_n=\left\{\begin{array}{cc} \frac{1}{2}\left(1-\left(\frac{n}{p}\right)\right), & \gcd(n,p)=1,\\
             0, & n\equiv 0\bmod p,
            \end{array}\right.
\end{equation}
where 
$$\left(\frac{n}{p}\right)=\left\{\begin{array}{cc} 1,& n \mbox{ a quadratic residue modulo }p,\\
-1, & n\mbox{ a quadratic non-residue modulo }p,\\
0, & n\equiv 0\bmod p,\end{array}\right.$$ 
is the {\em Legendre symbol}.
Since there are $(p-1)/2$ quadratic residues and $(p-1)/2$ quadratic non-residues ($0$ is neither a residue nor a non-residue), we obviously have 
$$B({\cal L}_p,p)=1.$$
By the Burgess bound, see for example \cite[$(12.58)$]{iwko}, we have
$$B({\cal L}_p,N)=O\left(N^{1-\frac{1}{r}}p^{\frac{r+1}{4r^2}}(\log p)^{\frac{1}{r}}\right)$$
for any $r\ge 1$. In particular, we have ($r=1$)
$$B({\cal L}_p,N)=O\left(p^{\frac{1}{2}}\log p\right)$$
and 
$$B({\cal L}_p,N)=o(N)\quad \mbox{for $N\ge p^{\frac{1}{4}+o(1)}$}.$$

\subsection{Two-prime generator}

For two odd primes $p$ and $q$ with, say, $p<q<2p$ the {\em two-prime generator} ${\cal W}=(w_n)_{n=0}^{\infty}$ of period $pq$ 
satisfies
$$w_n=\frac{1}{2}\left(1-\left(\frac{n}{p}\right)\left(\frac{n}{q}\right)\right),\quad \gcd(n,pq)=1.$$
For any choice of $w_n$ with $\gcd(n,pq)>1$, by \cite[Lemma~4]{brwi} we have
$$B({\cal W},N)=O\left((pq)^{1/2}\log(pq)\right),\quad 1\le N\le pq,$$
and again by the Burgess bound
$$B({\cal W},N)=o(N)\quad \mbox{for }N\ge (pq)^{1/4+o(1)}.$$

\subsection{Thue-Morse sequence (along squares)}

The {\em Thue-Morse sequence} ${\cal T}=(t_n)_{n=0}^\infty$ over $\F_2$ is defined by
\begin{equation}\label{tmdef}t_n=\left\{ \begin{array}{cl} t_{n/2}, & \mbox{$n$ even},\\ t_{(n-1)/2}+1, & \mbox{$n$ odd},
              \end{array}\right.\quad n=1,2,\ldots
\end{equation}
with initial value $t_0=0$.
Since $t_{2n}\neq t_{2n+1}$ we have 
$$B({\cal T},N)\le 1,\quad N=1,2,\ldots$$
This already points to some undesirable structure of the Thue-Morse sequence. Further weaknesses of this sequence are mentioned below.
We will also see that the Thue-Morse sequence has some desirable features such as a large linear complexity.

Certain sub-sequences, such as the sub-sequence of the Thue-Morse sequence along squares, may keep
the good properties of the original sequence but avoid the bad ones.
For the sub-sequence of the {\em Thue-Morse sequence along squares} ${\cal Q}=(t_{n^2})_{n=0}^\infty$
we have
$$B({\cal Q},N)=o(N)$$
by Mauduit and Rivat \cite[Th\'eor\`eme~1]{mari}.

\subsection{Omega sequence (for integers and polynomials)}

Let $n=p_1^{a_1}p_2^{a_2}\cdots p_r^{a_r}$ be the prime factorization of a positive integer $n$. The $\Omega$ function is defined by
$$\Omega(n)=a_1+a_2+\ldots+a_r.$$
We consider the sequence ${\cal O}=(o_n)_{n=0}^\infty$ with
$$o_0=0,\quad  o_n=\Omega(n)\bmod 2, \quad n=1,2,\ldots$$
We have 
$$B({\cal O},N)=o(N)$$
and the Riemann hypothesis is equivalent to, see Humphries~\cite{hu},
$$B({\cal O},N)=O\left(N^{\frac{1}{2}+\varepsilon}\right)\quad \mbox{for any $\varepsilon>0$}.$$

Similarly, for a polynomial $F(X)$ over the finite field~$\F_p$ of prime order $p$, $\Omega_p(F)$ denotes the total number of 
irreducible factors over $\F_p$ of $F(X)$. For fixed degree $d\ge 3$ we order the monic polynomials of degree $d$,
$$F_n(X)=X^d+n_{d-1}X^{d-1}+\ldots+n_1X+n_0,$$
where 
$$n=n_0+n_1p+\ldots+n_{d-1}p^{d-1}\quad \mbox{with } 0\le n_0,n_1,\ldots,n_{d-1}<p,$$ and define the sequence
${\cal P}={\cal P}_{d,p}=(p_n)_{n=0}^{p^d-1}$ of length $p^d$ by
$$p_n=\Omega_p(F_n)\bmod 2,\quad n=0,1,\ldots,p^d-1.$$
Carlitz \cite{ca} proved for $N=p^d$
$$B({\cal P},p^d)=p^{\left\lfloor\frac{d+1}{2} \right\rfloor}$$
and for $p>2$ and $N<p^d$ we have, see \cite{mewi16},
$$B({\cal P},N)=O\left(\frac{dN}{p^{1/2}}\right) \quad \mbox{for } N\ge p^2\log p,$$
that is, for fixed $d\ge 3$ and $p\rightarrow \infty$ we get
$$B({\cal P},N)=o(N),\quad N\ge p^2\log p.$$

A closely related measure of pseudorandomness, the {\em well-distribution measure}, that is, roughly speaking, the balance of the sequence along arithmetic progressions, was studied in the series of papers
\cite{al,camasa,masa,masa98,mewi16,mewi22,risa}.

\section{Linear complexity}
\subsection{Definition}
 The {\em $N$th linear complexity} $L({\cal S},N)$ of a binary sequence ${\cal S}$ 
 is the smallest positive integer $L$ such that there are constants 
$c_0,\ldots,c_{L-1}\in \F_2$ with
$$s_{n+L}=c_{L-1}s_{n+L-1}+\ldots+c_0s_n,~0\le n<N-L.$$
 The {\em linear complexity} $L({\cal S})$ of ${\cal S}$ is 
 $$L({\cal S})=\sup_{N\ge 1} L({\cal S},N).$$
 In particular, for a $T$-periodic sequence ${\cal S}_T$ we have
 $$L({\cal S}_T)\le T$$
 and $L({\cal S})<\infty$ if and only if ${\cal S}$ is ultimately periodic.

A sequence of small linear complexity is predictable and thus unsuitable in cryptography. However, the converse is not true.
There are many predictable sequences of very large linear complexity, for example periodic sequences containing only a single one in a period, and, in addition, finer quality measures have to be studied.

\subsection{Expected value}

Let denote by $A(N,L)$ the number of $(s_0,\ldots,s_{N-1})\in \F_2^N$ which are the initial values of a sequence ${\cal S}$ with $L({\cal S},N)=L$. The expected value 
$$E_N=\frac{1}{2^N}\sum_{L=0}^N A(N,L)L$$
was analyzed in Gustavson \cite{gu}.

\begin{theorem}
  The expected value 
  of $L({\cal S},N)$ is
  $$E_N=\frac{N}{2}+O(1).$$
\end{theorem}

Niederreiter~\cite{ni} showed that the $N$th linear complexity of a random sequence follows closely but
irregularly the $N/2$-line and deviations from $N/2$ of the order of magnitude $\log N$ must appear
for infinitely many $N$. 

From a computational point of view to avoid an attack via the {\em Berlekamp-Massey algorithm}, see \cite{ma69}, say $L({\cal S},N)\ge N^{o(1)}$ would be good enough.

For periodic sequences the expected value of the linear complexity depends on the period \cite{meni}.
For example, if the period $T$ is a prime and $2$ is a primitive root modulo $p$, then any non-constant sequence ${\cal S}$ of period $T$ is of linear complexity $T$ or $T-1$ and the expected value is very close to $T$, see \cite{cudire}.

\subsection{Legendre sequence and two-prime generator}

The linear complexity of the Legendre sequence ${\cal L}_p$ defined by $(\ref{legdef})$ was determined by Turyn \cite{tu}, see also \cite{dihesh}.

\begin{theorem}
  For a prime $p>2$ the linear complexity $L({\cal L}_p)$ of the $p$-periodic Legendre sequence~${\cal L}_p$ is
  $$L({\cal L}_p)=\left\{\begin{array}{cl}
   (p-1)/2, & p\equiv 1 \bmod 8,\\
   p, & p\equiv 3\bmod 8,\\
   p-1, & p\equiv -3 \bmod 8,\\
   (p+1)/2, & p\equiv -1\bmod 8.
   \end{array}\right.$$
\end{theorem}

For the $N$th linear complexity we have the following bound due to Chen et al.\ \cite{chgogoti}.
\begin{theorem}\label{leglin}
  $$L({\cal L}_p,N)\ge \frac{\min\{N,p\}}{p^{1/2}}\quad \mbox{for }N=0,1,\ldots$$
\end{theorem}
It would be important to improve this lower bound getting closer to the conjectured lower bound $N/2+o(N)$.

For the two-prime generator we get by \cite[Theorem~8.2.9]{cudire} 
$$L({\cal W})\ge \frac{(p-1)(q-1)}{2}$$
and by \cite{brwi}
$$L({\cal W},N)\ge \frac{\min\{N,pq\}}{(pq)^{1/2}}.$$

\subsection{Thue-Morse sequence (along squares)}

The following two results are due to \cite{mewi18} and \cite{suwi2}.
\begin{theorem}
For the $N$th linear complexity of the  Thue-Morse sequence ${\cal T}$ we have
$$L({\cal T},N)=2\left\lfloor \frac{N+2}{4}\right\rfloor,\quad N=1,2,\ldots$$
\end{theorem}

\begin{theorem}
For the $N$th linear complexity of the Thue-Morse sequence ${\cal Q}$ along squares
we have
$$L({\cal Q},N)\ge \left(\frac{2N}{5}\right)^{1/2}\quad \mbox{for }N\ge 21.$$
\end{theorem}

Note that the deviation of the $N$th linear complexity of the Thue-Morse sequence from $N/2$ is $O(1)$ which is too regular.
For the Thue-Morse sequence along squares we conjecture the desirable $L({\cal Q},N)=\frac{N}{2}+o(N)$.

\subsection{Omega sequences}

Up to our knowledge there is no lower bound on 
$L({\cal O},N)$ in the literature. However, our numerical data leads to the following conjecture.
\begin{conjecture}
$$L({\cal O},N)=\frac{N}{2}+O(1).$$
\end{conjecture}
If this conjecture is true, then the integer omega sequence can be distinguished from a random sequence by the deviation of the $N$th linear complexity from $N/2$.

For the polynomial omega sequence, combining \cite[Corollary~4]{chgogoti} and \cite{mewi16} we get
a lower bound on $L({\cal P},N)$
of order of magnitude $\min\{N,p^d\}^{1/2}p^{\frac{1}{4}-\frac{d}{2}}$.

\subsection{Balance and linear complexity}

Balance and linear complexity are independent measures of pseudorandomness in the following sense:
\begin{enumerate}
\item Both measures detect the non-randomness of constant sequences.
\item The non-randomness of the sequence $s_n=0$, $n=0,1,\ldots,N-2$, $s_{N-1}=1$, is detected by the balance but not by the $N$th linear complexity.
\item The balance of the Thue-Morse sequence is too small but its $N$th linear complexity is large enough.
\item We have seen several examples, for example the Legendre sequence, with both a high $N$th linear complexity and a desirable balance. 
\end{enumerate}

\section{Correlation measure}

\subsection{Definition and expected value}

 The {\em $N$th correlation measure of order~$k$} of ${\cal S}$ introduced by Mauduit and S\'ark\"ozy \cite{masa} is 
$$C_k({\cal S},N)=\max_{M,D}\left|\sum^{M-1}_{n=0}(-1)^{s_{n+d_1}}\cdots (-1)^{s_{n+d_k}}\right|,\quad k\ge 1,$$ 
where the maximum is taken over all $D=(d_1,d_2,\ldots,d_k)$ with integers satisfying
$0\le d_1<d_2<\cdots<d_k$ and $1\le M\le N-d_k$.

The correlation measure of order $k$ provides information about the independence of parts of the sequence and their shifts.
For a random sequence this similarity and thus the correlation measure of order $k$ is expected to be small.
More precisely, by \cite{al} we have the following result.

\begin{theorem}
For any $\varepsilon>0$ there exist an $N_0$ such that for all $N\ge N_0$
we have for a randomly chosen sequence ${\cal S}$
\begin{equation}\label{correxp}\frac{2}{5}\sqrt{N\log{N\choose k}}<C_k({\cal S},N)<\frac{7}{4}\sqrt{N\log{N\choose k}}
\end{equation}
with probability at least $1-\varepsilon$.
\end{theorem}
Hence, $C_k({\cal S},N)$ should be up to some logarithmic factor of order of magnitude $\sqrt{N}$ or at least $o(N)$.

\subsection{Correlation measure and linear complexity}

The following lower bound on the linear complexity profile in terms of the correlation measure was proved in \cite{brwi06}.

\begin{theorem}
Let ${\cal S}$ be a $T$-periodic binary sequence. For $2\le N\le T$
we have
$$ L({\cal S},N) \ge N-\max_{1\le k \le L({\cal S},N)+1}C_k(s_n,T). $$
\end{theorem}
For a recent improvement which saves typically a factor $\log N$ see Chen et al.\ \cite[Corollary 4]{chgogoti}.
For example, combining this relation between linear complexity and correlation measure with the bound on the correlation measure in Theorem~\ref{legcorr} below we immediately get the lower bound on the linear complexity of Theorem~\ref{leglin} above. In this sense we may say that the correlation measure of order $k$ is a finer measure of pseudorandomness than the linear complexity. 
However, from an algorithmic point of view the $N$th correlation measure of order~$k$ is much more difficult to analyze than the $N$th linear complexity. Still, for some special number-theoretic sequences such as the Legendre sequence one can estimate it theoretically.

\subsection{Legendre sequence and two-prime generator}

Although almost all sequences satisfy $(\ref{correxp})$, it is difficult to find concrete examples.
Roughly speaking, if you can describe a sequence, it does not behave like a randomly chosen sequence anymore.
However, for fixed $k$ and sufficiently large $p$, the correlation measure of order~$k$ of the Legendre sequence 
essentially behaves like the one for a randomly chosen sequence up to logarithmic terms, see \cite{masa}.

\begin{theorem}\label{legcorr}
  The correlation measure of order $k$ of the Legendre sequence satisfies
  $$C_k({\cal L}_p,N)=O(kp^{1/2}\log p),\quad 1\le N\le p.$$
\end{theorem}

The situation is different for the two-prime generator.
On the one hand, by \cite{risa} we still have 
$$C_2({\cal T},N)=O((pq)^{3/4}),\quad 1\le N\le pq.$$
On the other hand, taking the lags 
$$d_0=0, d_1=p, d_2=q, d_3=p+q$$ 
we get
$$C_4({\cal T},N)= N+O(N^{1/2}),\quad 1\le N\le pq,$$
showing that the two-prime generator is not a good candidate for cryptography.

\subsection{Thue-Morse sequence (along squares)}

By \cite{masa98} we have
$$C_2({\cal T},N)>\frac{N}{12},\quad N\ge 5.$$
We believe that this feature of non-randomness is destroyed by taking the sub-sequence along squares.
\begin{conjecture}
  $$C_k({\cal Q},N)=o(N)\quad \mbox{for }k=2,3,\ldots$$
\end{conjecture}
If we assume that the lags are bounded by a constant $C$, that is, $d_k\le C$ and $N\rightarrow \infty$, the analog of the correlation measure of order $k$ with bounded lags is $o(N)$ by
\cite{drmari}.

\section{Omega sequences}

The following is essentially Chowla's conjecture, see \cite{ch}.
\begin{conjecture}
  $$C_k({\cal O},N)=o(N).$$
\end{conjecture}
For recent progress on Chowla's conjecture see Tao and Ter\"av\"ainen \cite{tate1,tate2} and references therein. The correlation measure of order $k$ of a modified omega sequence was studied by Cassaigne et al.\ \cite{camasa}.

The Chowla conjecture for polynomials was settled by Carmon and Rudnick \cite{caru} for $p>2$ and Carmon \cite{car} for $p=2$ in the case that $p$ is fixed and the degree $d$ goes to infinity.
In particular, we have the following bound, see \cite[Theorem 3]{mewi16}. 
\begin{theorem}\label{Pdcorr}
 $$C_k({\cal P}_d,p^d)=O(k^2dp^{d-1/2}\log p).$$
\end{theorem}
However, nothing is known for fixed $p$ and $d\rightarrow \infty$. For polynomials over finite fields $\mathbb{F}_{p^r}$ with $r\ge 3$ there has been a recent breakthrough by Sawin and Shusterman \cite{sash}. 
However, it seems that the case $r=1$ is out of reach.

\section{Maximum-order complexity}

\subsection{Definition, expected value and relation to other measures}
The {\em $N$th maximum order complexity}  $M({\cal S},N)$ is the smallest positive integer~$M$ with 
 $$s_{n+M}=f(s_{n+M-1},\ldots,s_n),\quad 0\le n\le N-M-1,$$
 for some mapping $f:\F_2^M \mapsto \F_2$.
 
 The maximum order complexity was introduced by Jansen in \cite[Chapter 3]{ja1}, see
also \cite{ja2}. The typical value for the $N$th maximum order complexity is of order of magnitude $\log N$, see \cite{ja1,ja2}.

Obviously, we have 
$$M({\cal S},N)\le L({\cal S},N)$$
and we may consider the maximum-order complexity a finer measure than the linear complexity. However, from an algorithmic point of view the linear complexity can be much easier determined via the Berlekamp-Massey algorithm than the maximum-order complexity. 
An algorithm for calculating the maximum order complexity
profile of linear time and memory was presented by Jansen \cite{ja1,ja2} using the graph algorithm introduced by Blumer et al.~\cite{bl}.
 
 Although a large $M({\cal S},N)$ is desired it should not be too large since otherwise the correlation measure of order $2$
 is large, see \cite[(5.6)]{mewi22}. Combining this inequality with \cite[Theorem 5]{chgogoti} we get:
 \begin{theorem}\label{maxrel}
   We have 
   $$C_2({\cal S},N)\ge \max\left\{M({\cal S},N)-1,N+1-2^{M({\cal S},N)}\right\}.$$
 \end{theorem}
 
 In Subsection~\ref{motm} we study the {\em Thue-Morse sequence} ${\cal T}=(t_n)_{n=0}^\infty$
 defined by $(\ref{tmdef})$.
 It turns out that $M({\cal T},N)$ is of order of magnitude $N$. However, this implies that the correlation measure $C_2({\cal T},N)$ of order $2$
 is also of order of magnitude~$N$ and concerning this measure the Thue-Morse sequence does not behave like a random sequence.
 
 However, for the {\em Thue-Morse sequence along squares} ${\cal Q}=(t_{n^2})_{n=0}^\infty$
 we mention that $M({\cal Q},N)$ is at least of order of magnitude $N^{1/2}$. 
 
 We can also define and study
 $$M({\cal S})=\sup_{N\ge 1}M({\cal S},N),$$
 see for example \cite{liko}.

 \subsection{Legendre sequence}
 
 Combining Theorem~\ref{maxrel} and Theorem~\ref{legcorr} we get the following bound on the maximum-order complexity of the Legendre sequence.
 \begin{corollary}
  For $1\le N\le p$ we have 
  $$M({\cal L}_p,N)=O(p^{1/2}\log p)$$
  and
  $$M({\cal L}_p,N)\ge \frac{\log N-\frac{1}{2}\log p+O(\log\log p)}{\log 2}.$$ 
 \end{corollary}

It is not difficult to obtain a similar bound for the two-prime generator. However, because of its large correlation measure of order $4$ there is no need of further studies of this sequence.

\subsection{Thue-Morse sequence (along squares)}\label{motm}

The following result is due to \cite{suwi1}.
 \begin{theorem}
  For $N\ge 4$,
  the $N$th maximum order complexity of the Thue-Morse sequence ${\mathcal T}$
  satisfies\vadjust{\vskip-1.5ex}
\begin{gather*}
  M({\mathcal T},N)=2^\ell+1,
\\[-1ex]
\noalign{\leftline{where}}
  \ell=\left\lceil \frac{\log (N/5)}{\log 2}\right\rceil.
\end{gather*}
 \end{theorem}
 It is easy to see that
\begin{gather*}
 \frac{N}{5}+1\le M({\mathcal T},N)\le 2\frac{N-1}{5}+1
 \quad\text{for}\;\; N\ge 4
\\[-1ex]
 \noalign{\leftline{and\vadjust{\vskip-0.5ex}}}
 M({\mathcal T},1)=0,
\quad
 M({\mathcal T},2)=M({\mathcal T},3)=1.
\end{gather*}

For the Thue-Morse sequence along squares see \cite{suwi2}.
 \begin{theorem}
    $$ M({\cal Q},N)\geq \sqrt{\frac{2N}{5}},\quad N\ge 21.$$
\end{theorem}
For an extension to sub-sequences of the Thue-Morse sequence along polynomial values see \cite{po}.

\subsection{Omega sequences}

Our numerical data leads to the following conjecture.
\begin{conjecture}
 $M({\cal O},N)$ is of order of magnitude $\log N$.
\end{conjecture}

For the polynomial analog we get from Theorems~\ref{Pdcorr} and~\ref{maxrel} the following result.
\begin{corollary}
$$M({\cal P}_d,p^d)=O(dp^{d-1/2}\log p)$$
and
$$M({\cal P}_d,p^d)\ge \frac{(d-1/2)\log p}{\log 2}+o(\log p).$$
\end{corollary}

\section{Expansion complexity}

\subsection{Definition and Thue-Morse sequence}

 Let 
 $$G(x)=\sum_{n=0}^\infty s_nx^n$$ 
 be the {\em generating function} of the sequence ${\cal S}$ with 
 $(s_0,\ldots,s_{N-1})\ne (0,\ldots,0)$. The smallest degree $E({\cal S},N)$ of a polynomial $h(x,y)\ne 0$ with
 $$h(x,G(x))\equiv 0\bmod x^N$$
 is called {\em $N$th expansion complexity of ${\cal S}$}.\\
 For $(s_0,\ldots,s_{N-1})=(0,\ldots,0)$ we define $E({\cal S},N)=0$.\\
 The {\em expansion complexity} $E({\cal S})$ of ${\cal S}$ is
 $$E({\cal S})=\sup_{N\ge 1} E({\cal S},N).$$
 The expansion complexity was introduced by Diem \cite{di} and we have
 $$E({\cal S},N)<\sqrt{2N}$$
 by \cite[Theorem 4]{gomeni}.
 By Christol's theorem \cite{chr} {\em automatic sequences} are characterized by $E({\cal S})<\infty$.

From the well-known equation
$$(1+x)^3G(x)^2+(1+x)^2G(x)+x=0$$
we immediately get the following bound.
\begin{corollary}
 For $N=1,2,\ldots$, the $N$th expansion complexity $E({\cal T},N)$ of the Thue-Morse sequence is at most~$5$. 
\end{corollary}

The expansion complexity is another measure for the predictability of a sequence.
Despite of its very large $N$th linear complexity,
the Thue-Morse sequence is very predictable because of its extremely small $N$th expansion complexity.
Hence, the $N$th expansion complexity can be substantially smaller than the $N$th linear complexity. However, we will see in the next section that in the periodic case expansion complexity and linear complexity are essentially the same.

The expected value of $E({\cal S},N)$ of a random sequence is of order of magnitude $N^{1/2}$, see \cite[Theorem 2]{gome}.

Diem showed \cite{di} that if a sequence has small expansion complexity, then long parts of
such sequences can be computed efficiently from short ones. An algorithm based on Gr\"obner basis is given in \cite{gome}.

\subsection{Expansion complexity and linear complexity}
The following results are from \cite{meniwi}.

\begin{theorem}
 Let ${\cal S}$ be a (purely) periodic sequence. Then we have
 $$E({\cal S})=L({\cal S})+1.$$
\end{theorem}

In the aperiodic case we get the following.
\begin{theorem}
 $$E({\cal S},N)\le \min\{L({\cal S},N)+1,N+2-L({\cal S},N)\}.$$
\end{theorem}
This result and the results on the Thue-Morse sequence imply that the $N$th expansion complexity is a strictly finer measure than the linear complexity, more precisely, than the deviation of $L({\cal S},N)$ from the expected value $\frac{N}{2}$.

\subsection{Thue-Morse along squares, Legendre sequence and omega sequences}

Our numerical data leads to the conjecture that all four sequences have expansion complexity of order of magnitude~$N^{1/2}$ with the restriction $N\le p^2$ for the Legendre sequence.

\section{$2$-adic complexity}

Besides the linear complexity, that is the length of a shortest linear feed shift register which generates the sequence, the $2$-adic complexity has been studied, which is closely related to the length of a shortest 
feedback with carry shift registers which generates the sequence and was introduced by Goresky and Klapper, see \cite{gokl} and references therein. Although the theory of $2$-adic complexity has been very well-developed for the periodic case, almost nothing is known for the aperiodic case. 

More precisely, the {\em $2$-adic complexity} $C({\cal S})$ of a $T$-periodic sequence is 
$$C({\cal S})=\frac{\log\left(\frac{2^T-1}{\gcd(2^T-1,S(2)}\right)}{\log 2},$$
where 
$$S(X)=\sum_{n=0}^{T-1}s_nX^n.$$
The expected value of the $2$-adic complexity of $T$-periodic sequences is $T-O(\log T)$, see \cite[Corollary~18.2.2]{gokl}.

Since the linear complexity satisfies
$$L({\cal S})=T-\deg(\gcd(S(X),X^T-1))$$
it is easy to see that linear complexity and $2$-adic complexity complement each other. 
For example, let $2^T-1$ be a Mersenne prime. Then any non-constant sequence has maximum $2$-adic complexity. However, $X^T-1$ may still have a nontrivial divisor $S(X)$ of large degree and the linear complexity can be small.
Conversely, if $T$ is a prime and $1+X+\ldots+X^{T-1}$ is irreducible, that is, $2$ is a primitive root modulo $T$, then any non-constant sequence has maximal linear complexity. However, $2^T-1$ may have a large nontrivial divisor and the $2$-adic complexity can be small. 

Moreover, an {\em $m$-sequence} of period $2^r-1$ has linear complexity only $r$ but maximal $2$-adic complexity \cite{XQL}. Conversely, {\em $\ell$-sequences} have minimal $2$-adic complexity but can have very large linear complexity 
\cite{QX}.

The Legendre sequence has maximal $2$-adic complexity and the $2$-adic complexity of the two-prime generator is very large as well if $p$ and $q$ are essentially of the same size, see \cite{howi,hhu,XQL}.  

It would be very important to study also the aperiodic case, in particular, to get results for the Thue-Morse sequence along squares and the omega sequences. 
More precisely, the $N$th $2$-adic complexity $C({\cal S},N)$ of ${\cal S}$ is the binary logarithm of
$$\min\{\max\{|f|,|g|\}: f,g\in \mathbb{Z}, g\mbox{ odd, }gS(2)\equiv f\bmod 2^N\},$$
where 
$$S(2)=\sum_{n=0}^{N-1}s_n2^n.$$
The question about the expected value is open, see \cite[Section~5.1]{kl}.
Our numerical data obtained using the {\em rational approximation algorithm}, see \cite[Chapter 17]{gokl}, leads, for example, to the following conjecture for the Legendre sequence.
\begin{conjecture}
  $$C({\cal L}_p,N)=\frac{\min\{N,2p\}}{2}+O(1).$$
\end{conjecture}
Similar conjectures can be stated for the Thue-Morse sequence along squares and the omega-sequences.

Another very promising balanced number-theoretic sequence ${\cal X}$ defined by 
$$x_n=(g^n \bmod p)\bmod 2,\quad n=0,1,\ldots$$ 
for some $g\in \F_p^*$ has been introduced very recently in \cite{pa}. 
It is natural to ask for lower bounds on linear complexity, maximum-order complexity etc.\ for this sequence as well. Note that in the case that $g=2$ is a primitive root modulo $p$ we get an $\ell$-sequence by \cite[Theorem 4.5.2]{gokl}.

There is also a $2$-adic analog of the correlation measure of order $2$ (aperiodic autocorrelation) called {\em arithmetic autocorrelation} which can be estimated in terms of the correlation measure, see \cite{homewi}.
In particular, the arithmetic autocorrelation of the Legendre sequence was estimated in \cite{howi2}. 

\section*{Acknowledgment}
The author wishes to thank Zhixiong Chen and L\'aszl\'o M\'erai for useful comments.

\end{document}